\documentclass[8pt,letterpaper,twocolumn]{amsart}
\usepackage{ragged2e} 
\usepackage{amsmath, amssymb, graphicx, url, amsthm,  enumitem, textcomp, gensymb, mathtools, braket}
\usepackage{mdframed}
\usepackage{indentfirst}
\usepackage{units}
\usepackage{wrapfig}
\usepackage[dvipsnames]{xcolor}

\usepackage{tikz}
\usepackage{caption}
\usetikzlibrary{backgrounds,shapes.geometric,arrows}
\usepackage{tcolorbox}
\usepackage[T1]{fontenc}
\usepackage{titlesec}
\usepackage{enumitem}

\usepackage[dvips,centering,width=15.1cm,height=21.3cm,
paperheight=23.8cm,paperwidth=17cm
]{geometry}

\titleformat{\section}
{\color{purple}\bfseries\scshape\centering}
{\thesection.}{0.5em}{}

\titleformat{\subsection}
{\bfseries\centering}
{\thesection.}{0.5em}{}

\newcommand{\R}{\mathbb{R}}
\newcommand{\C}{\mathbb{C}}

\newcommand{\bfS}{\mathbf{S}}

\newcommand{\DD}{\mathbb{D}}
\newcommand{\HH}{\mathbb{H}}

\let\ol\overline

\usepackage{etoolbox}
\AtBeginEnvironment{proof}{\setlength{\parindent}{0.4cm}}
\AtBeginEnvironment{itemize}{\setlength{\parindent}{0.2cm}}
\AtBeginEnvironment{enumerate}{\setlength{\parindent}{0.2cm}}

\theoremstyle{definition}
\newtheorem{definition}{Definition}
\newtheorem{remark}[definition]{Remark}

\newtheorem{observation}[definition]{Observation}

\theoremstyle{plain}
\newtheorem{theorem}[definition]{Theorem}

\newtheorem{proposition}[definition]{Proposition}
\newtheorem*{proposition-nonum}{Proposition}
\newtheorem*{lemma-nonum}{Lemma}
\newtheorem*{theorem-nonum}{Theorem}
\newtheorem*{example-nonum}{Example}
\newtheorem*{definition-nonum}{Definition}

\newcommand{\lr}[1]{\langle #1 \rangle}

\usepackage{hyperref}

\hypersetup{
  colorlinks   = true, 
  urlcolor     = blue, 
  linkcolor    = blue,
  citecolor   = red 
}

\title{\Huge Ptolemy's equation \\[5pt] and kin}
\author{\Large Katie Waddle}
\address{Department of Mathematics, University of Michigan, Ann Arbor, MI 48109, USA}
\email{waddle@umich.edu}
\thanks
{\emph{2020 Mathematics Subject Classification}
Primary 51K99,  
Secondary 51M25, 
97G40. 
}

\thanks{Partially supported by National Science Foundation grant DMS-2348501.}

\keywords{Distance geometry, 3-term relations, cross-ratio.}

\date{\today}

\begin{document}

\pagestyle{plain}

\maketitle

\newpage

Three-term relations of the form 
$AB+CD=EF$ 
arise in multiple mathematical contexts, including the Ptolemy equation for a cyclic quadrilateral, Casey's theorem on bitangents, Penner's relation for lambda lengths, and  Pl{\"u}cker's identity for the maximal minors of a~$2\times 4$ matrix.  
In this note, we explain how these different occurrences of the 3-term relation can be directly obtained from each other.

\section*{Four incarnations of the 3-term relation}\label{sec: four occurrences}

\begin{figure}[htbp!]
    \centering
\begin{tikzpicture}[scale=2,anode/.style={scale=1},rotate=90]
\draw (0,0) circle (1);
    \coordinate (a1) at ({1/sqrt(2)},{1/sqrt(2)});
    \node at (a1) [left,anode] {$ A_1$};
            
    \coordinate (a2) at (0,1);
    \node at (a2) [left,anode] {$ A_2$};
           
    \coordinate (a3) at ({-sqrt(2)/sqrt(3)},{-1/sqrt(3)}) {};
    \node at (a3) [right,anode] {$ A_3$};
             
    \coordinate (a4) at ({1/sqrt(2)},{-1/sqrt(2)});
            
    \node at (a4) [right,anode] {$ A_4$};
    \draw[cyan, very thick] (a1) to node[right, outer sep=-1pt] {$d_{12}$} (a2) to node[left] {$d_{23}$} (a3) to node[right] {$d_{34}$} (a4) to node[above,inner sep=1pt] {$d_{14}$} (a1) to node[right] {$d_{13}$} (a3);
        
    \draw[cyan, very thick] (a2) to node[below, near end, inner sep=2pt] {$d_{24}$}  (a4);

    \filldraw (a1) circle (.02);
            \filldraw (a2) circle (.02);
            \filldraw (a3) circle (.02);
            \filldraw (a4) circle (.02);
\end{tikzpicture}
\vspace{-5pt}
\caption{Distances between four points on a circle}
\label{fig: ptolemy}
\end{figure}
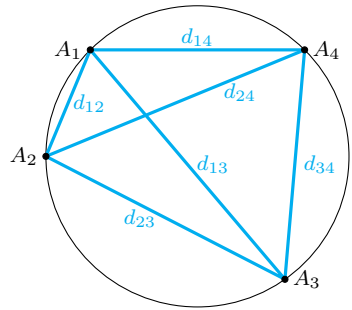

\textbf{Ptolemy}, a Greco-Roman astronomer and geometer, lived in the city of Alexandria in the 2nd century CE. His classical treatise \emph{Almagest} contains the following famous result:

\begin{tcolorbox}
\begin{theorem}[{\cite[Book~1, Sec.~10, Lem.~1.1]{ptolemyAlmagestIntroductionMathematics2014}}] \label{thm: ptolemy} Let $A_1,A_2,A_3,A_4$ be four points  arranged counterclockwise on a unit circle~$\bfS$ on the Euclidean plane. 
For $1\le i<j\le 4$, let~$d_{ij}$ denote the distance from~$A_i$ to~$A_j$, see Figure~\ref{fig: ptolemy}. Then 
\begin{equation} \label{eq: ptolemy}
d_{12}d_{34}+d_{14}d_{23}=d_{13}d_{24}.
\end{equation}
\end{theorem}

\end{tcolorbox}

\begin{figure}[htbp!]
    \centering
    \begin{tikzpicture}[scale=3,dnode/.style={scale=1},anode/.style={scale=1},rotate=90]
\draw (0,0) circle (1);
    \coordinate (a1) at ({1/sqrt(2)},{1/sqrt(2)});
   
    \node at (a1) [left,anode] {$ A_1$};
    \coordinate (a2) at (0,1);

    \node at (a2) [left,anode] {$ A_2$};
    \coordinate (a3) at ({-sqrt(2)/sqrt(3)},{-1/sqrt(3)}) {};
    \node at (a3) [right,anode] {$ A_3$};

    \coordinate (a4) at ({1/sqrt(2)},{-1/sqrt(2)});
    \node at (a4) [right,anode] {$ A_4$};

\draw ({3/(4*sqrt(2))},{3/(4*sqrt(2))}) circle (.25);

\draw (0,.8) circle (.2);

\draw ({-sqrt(2)/sqrt(3)+sqrt(2)/(3*sqrt(3))},{-1/sqrt(3)+1/(3*sqrt(3))}) circle ({1/3});

\draw ({5/(8*sqrt(2)},{-5/(8*sqrt(2))}) circle ({3/8});

\draw[purple, very thick] (-.55,-.72) --  node[right, inner sep=1pt, outer sep=2pt,dnode] {$t_{34}$}(.43,-.815);

\draw[purple, very thick] (.82,-.45) --  node[above, inner sep=2pt,dnode] {$t_{14}$}(.78,.56);

\draw[purple, very thick] (.055,.999) --  node[right, inner sep = 2pt,pos=.45,dnode] {$t_{12}$}(.65,.75);

\draw[purple, very thick] (-.85,-.25) --  node[left,inner sep = 4pt,dnode] {$t_{23}$}(-.17,.9);

\draw[purple, very thick] (-.20,.78) -- node[below, inner sep=4pt,dnode,pos=.7] {$ t_{24}$}
 (.07,-.5);

 \draw[purple, very thick] (-.75,-.125) -- node[right,pos=.7,inner sep=4pt,dnode] {$t_{13}$}(.4,.745);

 \filldraw (a1) circle (.02);\filldraw (a2) circle (.02);\filldraw (a3) circle (.02); \filldraw (a4) circle (.02);

 \node at (.6,.4) {$H_1$};
\node at (0,.7) {$H_2$};
\node at (-.4,-.2) {$H_3$};
\node at (.4,-.2) {$H_4$};
\end{tikzpicture}
\caption{Bitangent distances between four circles tangent to a fifth circle}
\label{fig: casey}
\end{figure}
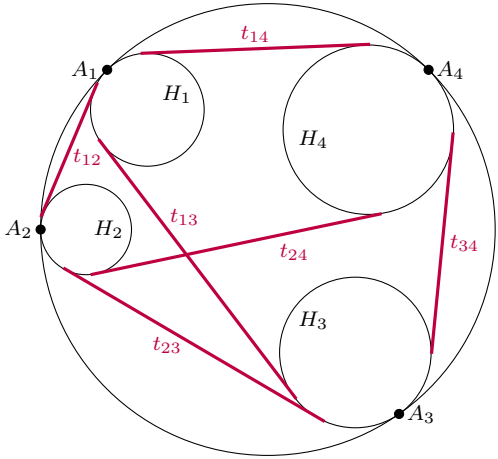

In 1820s Japan, \textbf{Shiraishi Nagatada} (also known as Shiraishi Ch\=och\=u) depicted a remarkable generalization of Ptolemy's theorem without providing a proof (see \cite{fukagawaSacredMathematics2008}). 
In 1860s Ireland, \textbf{John Casey}, at the time a high school teacher, independently published this result:
\begin{tcolorbox}
\RaggedRight
\begin{theorem}[{\cite{caseyEquationsProperties11866}}]\label{thm: casey}
Let $H_1,H_2,H_3,H_4$ be 
four circles contained in~$\bfS$ that are tangent to~$\bfS$ at the points~$A_1,A_2,A_3,A_4$, respectively,
see Figure~\ref{fig: casey}. 
For $1\le i<j\le 4$, let $t_{ij}$ denote the length of a common tangent to~$H_i$ and $H_j$ drawn so that both circles lie on the same side. 
Then
\begin{equation} \label{eq: casey}
t_{12}t_{34}+t_{14}t_{23}=t_{13}t_{24}.
\end{equation}
\end{theorem}

\end{tcolorbox}
\noindent In the degenerate case, when the four circles shrink to points on~$\bfS$, we recover Ptolemy's theorem.

\pagebreak[3]

\begin{figure}[htbp!]
\centering
    \begin{tikzpicture}[scale=3,dnode/.style={scale=1},rotate=90]
\draw (0,0) circle (1);
    \coordinate (a1) at ({1/sqrt(2)},{1/sqrt(2)});
    \filldraw (a1) circle (.02);
    \node at (a1) [left] {$ A_1$};
    \coordinate (a2) at (0,1);
        \filldraw (a2) circle (.02);

    \node at (a2) [left] {$ A_2$};
    \coordinate (a3) at ({-sqrt(2)/sqrt(3)},{-1/sqrt(3)}) {};
    \node at (a3) [right] {$ A_3$};
        \filldraw (a3) circle (.02);

    \coordinate (a4) at ({1/sqrt(2)},{-1/sqrt(2)});
    \node at (a4) [right] {$ A_4$};
        \filldraw (a4) circle (.02);

\draw ({3/(4*sqrt(2))},{3/(4*sqrt(2))}) circle (.25);

\draw (0,.8) circle (.2);

\draw ({-sqrt(2)/sqrt(3)+sqrt(2)/(3*sqrt(3))},{-1/sqrt(3)+1/(3*sqrt(3))}) circle ({1/3});

\draw ({5/(8*sqrt(2)},{-5/(8*sqrt(2))}) circle ({3/8});

\draw[black] (a1) to[out=225,in=270]  (a2);

\draw[very thick, brown] (.31,.65) to [out=156,in=-35] node[dnode,  left, pos=.1, inner sep=5pt] {$\delta_{12}$} (.18,.71);

\draw[black] (a2) to[out=270,in=31] (a3);

\draw[very thick, brown] (-.031,.6) to [out=261,in=55] node[below,dnode] {$\delta_{23}$}  (-.32,-.14);

\draw[black] (a3) to[out=31,in=135]   (a4);

\draw[very thick, brown] (-.21,-.37) to [out=5,in=175] node[right,dnode] {$\delta_{34}$} (.075,-.38);

\draw[black] (a4) to[out=135,in=225]  (a1);

\draw[very thick, brown] (.455,.29) to [out=254,in=95] node[above,dnode] {$\delta_{14}$} (.42,-.07);

\draw[black] (a1) to[out=225,in=31]  (a3);

\draw[very thick, brown] (.355,.355) to [out=225,in=40] node[right, pos=.6, inner sep=5pt,dnode]{$\delta_{13}$} (-.262,-.2);

\draw[black] (a2) to[out=270,in=135] (a4);

\draw[very thick, brown] (.025,.6) to [out=278,in=117] node[near start,above,inner sep=4pt,dnode] {$\delta_{24}$} (.248,-.118);

\node at (.8,.3) {$H_1$};
\node at (-.3,.7) {$H_2$};
\node at (-.8,0) {$H_3$};
\node at (.1,-.8) {$H_4$};

\end{tikzpicture}
\caption{Hyperbolic distances between horocycles}
\label{fig: penner}
\end{figure}
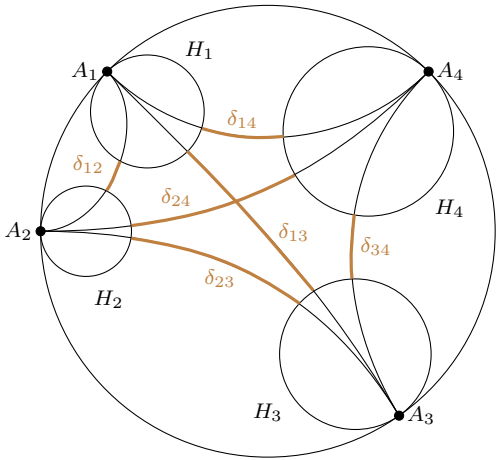

In the late 20th century, the American mathematician \textbf{Robert Penner} introduced the notion of \emph{lambda length}
\[
\lambda_{ij}=e^{\delta_{ij}/2}
\]
between two horocycles~$H_i$ and~$H_j$ in the hyperbolic plane; here~$\delta_{ij}$ is the signed hyperbolic distance between $H_i$ and~$H_j$. He established the following hyperbolic analogue of Ptolemy's theorem:
\begin{tcolorbox}
   \begin{theorem}[{\cite{penner_decorated_1987}}]\label{thm: penner}
Let $A_1,A_2,A_3,A_4$ be ideal points on the boundary of the hyperbolic plane. 
For $1\le i<j\le 4$, let $H_1,H_2,H_3,H_4$ be four horocycles centered at $A_1,A_2,A_3,A_4$, see Figure~\ref{fig: penner}. 
Then the lambda lengths $\lambda_{ij}$ satisfy
    \begin{equation}\label{eq: penner}
\lambda_{12}\lambda_{34}+\lambda_{23}\lambda_{14}=\lambda_{13}\lambda_{24}.
    \end{equation}
    \end{theorem}

\end{tcolorbox}

   Working at the same time as Casey, the German scientist \textbf{Julius Pl{\"u}cker} initiated the study of what we now call the \emph{Grassmannian} of 2-dimensional subspaces of a 4-dimensional vector space.  
  In modern notation, he discovered the following \emph{Pl\"ucker relation}:
  
\begin{tcolorbox}
  \begin{theorem}[{\cite[Add.~Note, Sec.~I, Part~3]{plscriptuscriptcker_new_1865}}] 
  \label{thm: plucker} 
  Consider a~$2\times 4$ matrix \[
  M=\left[\begin{matrix}
        x_1 & x_2 & x_3 & x_4\\
        y_1 & y_2 & y_3 & y_4
    \end{matrix}\right].\]
    For $1\le i<j\le 4$, denote \[{P_{ij}=\det\left[\begin{matrix}
      x_i & x_j\\
      y_i & y_j
  \end{matrix}\right]}.\]
  Then
  \begin{equation}\label{eq: plucker}
  P_{12}P_{34}+P_{14}P_{23}=P_{13}P_{24}.
  \end{equation}
Furthermore, any 6-tuple of complex numbers $(P_{ij})$ satisfying \eqref{eq: plucker} can be obtained in this way. 
\end{theorem}
\end{tcolorbox}
Many mathematicians have observed the striking similarity between the equations~\eqref{eq: ptolemy}-\eqref{eq: plucker}.  
See, e.g., A.~Felikson's survey~\cite{felikson_ptolemy_2023}, which also discusses connections to \emph{frieze patterns} and \emph{cluster algebras}. It is natural to wonder whether the equations~\eqref{eq: ptolemy}-\eqref{eq: plucker} can be directly deduced from each other.  We will show that this is indeed the case.

\vspace{.2cm}

\section*{Each equation implies the others}\label{sec: implications}

\vspace{.2cm}

We will employ the following  insight.

\begin{tcolorbox}
\RaggedRight
\begin{observation}
\label{obs: scaling} 
Consider two $6$-tuples of complex numbers:
$    (a_{12}, a_{13}, a_{14}, a_{23}, a_{24}, a_{34}),$ and $  (b_{12}, b_{13}, b_{14}, b_{23}, b_{24}, b_{34})$.
Suppose that there exist nonzero numbers~$q_1, q_2, q_3, q_4$ such that
\begin{equation}
\label{eq:ab-rescaling}
b_{ij} = q_i q_j a_{ij} \quad \text{(for~$1\le i<j\le 4$).}
\end{equation}
The numbers~$a_{ij}$ satisfy the 3-term relation
\begin{equation}
\label{eq:aij}
a_{12} a_{34} + a_{14} a_{23} = a_{13} a_{24}
\end{equation}
if and only if the numbers~$b_{ij}$ satisfy the 3-term relation
\begin{equation}
\label{eq:bij}
b_{12} b_{34} + b_{14} b_{23} = b_{13} b_{24} .
\end{equation}
\end{observation}
\end{tcolorbox}

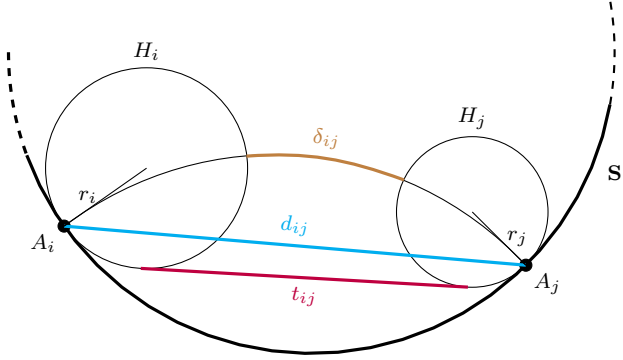
\begin{figure*}[htbp!]
\vspace{-5pt}
    \centering
\begin{tikzpicture}[scale=4]
    \draw [black,very thick,domain=200:350] plot ({cos(\x)}, {sin(\x)});
    \draw [black,very thick, dashed,domain=180:200] plot ({cos(\x)}, {sin(\x)});
    \draw [black,thick,dashed,domain=350:370] plot ({cos(\x)}, {sin(\x)});
    \node at (.95,-.4) [right] {$\bfS$};

    \draw ({-sqrt(8)/sqrt(27)},{-2/sqrt(27)}) circle ({1/3});
    \coordinate (A_1) at ({-sqrt(2)/sqrt(3)},{-1/sqrt(3)});
    \node at (A_1) [below left] {$A_i$};
    \filldraw (A_1) circle (.02);
    \node at (-{sqrt(8)/sqrt(27)},{-2/sqrt(27)+1/3}) [above] {$H_i$}; 

    \draw[] ({3/(sqrt(32))},{-3/(sqrt(32))}) circle (1/4);
    \coordinate (A_2) at ({1/sqrt(2)},{-1/sqrt(2)});
    \node at (A_2) [below right] {$A_j$};
            \filldraw (A_2) circle (.02);
    \node at ({3/(sqrt(32))},{-3/(sqrt(32))+1/4}) [above] {$H_j$};

    \coordinate (tan1) at (.5159405562764,-0.7799156252096);
    \coordinate (tan2) at (-0.5635170934365,-0.7176808985527);
     \draw[very thick, purple] (tan1) to node[below] {$t_{ij}$} (tan2);
    
    \draw[very thick, cyan] (A_1) to node[above] {$d_{ij}$} (A_2);

    \draw (A_1) to [out=37,in=133]  (A_2);
    \draw[very thick, brown] (-.215,-.345) to [out=4,in=160] node[above] {$\delta_{ij}$}(0.305,-0.424);

    \draw (A_1) to node[left] {$r_i$} ({-sqrt(8)/sqrt(27)},{-2/sqrt(27)});

    \draw (A_2) to node[right] {$r_j$} ({3/(sqrt(32))},{-3/(sqrt(32))});

\end{tikzpicture}    
\vspace{-2pt}
\caption{The measurements~$d_{ij}$,~$t_{ij}$, and~$\delta_{ij}$}
    \label{fig: distances}
\vspace{-2pt}
\end{figure*}

We next describe the measurements that will be cast in the roles of~$a_{ij},b_{ij}$ (see Figure~\ref{fig: distances}). 

Let~$\bfS\subseteq\R^2$ be the unit circle centered at the origin. Let~$A_1,A_2,A_3,A_4\in\bfS$ be four points on~$\bfS$ arranged in order counterclockwise.  For~${1\le i\le 4}$, let~${0\le \alpha_i\le\pi}$ be such that~$A_i=(\cos2\alpha_i,\sin2\alpha_i)$.  Let~$H_1, H_2, H_3, H_4$ be four circles tangent to $\bfS$ on the inside at $A_1, A_2, A_3, A_4$, as in Figure~\ref{fig: casey}. 
Let $r_i$ be the radius of~$H_i$. For simplicity we assume that~$H_i\cap H_j=\varnothing$ for~$1\le i,j\le 4$.  

If we change our perspective and regard~$\bfS$ as the Poin\-car{\'e} disk representing the hyperbolic plane,
then the~$H_i$ become horocycles with centers~$A_i$. Draw the hyperbolic geodesic $\gamma_{ij}$ from~$A_i$ to~$A_j$ and denote by~$\delta_{ij}$ the hyperbolic distance along~$\gamma_{ij}$ between the points~$H_i\cap\gamma_{ij}$ and~$H_j\cap \gamma_{ij}$. 

We now have the following measurements:
\begin{itemize}[leftmargin=.15in]
        \item~$d_{ij}=$  Euclidean distance between~$A_i$ and~$A_j$,
        \item~$t_{ij}=$ length of the exterior bitangent to $H_i$ and~$H_j$,
        \item~$\lambda_{ij}=e^{\delta_{ij}/2}=$  lambda length between~$H_i$ and~$H_j$,
        \item~$P_{ij}=\det\left[\begin{matrix}
    \cos\alpha_i & \cos\alpha_j\\
    \sin\alpha_i & \sin\alpha_j
\end{matrix}\right]$.
    \end{itemize}

It turns out that these measurements are related to each other by rescaling: 

\begin{tcolorbox}
\begin{proposition}\label{prop: scaling}
In the setting described above,
    \begin{align}
        t_{ij}&=\sqrt{1-r_i}\sqrt{1-r_j}d_{ij}\label{eq: d and t}\\
t_{ij}&=\lambda_{ij}\sqrt{2r_i}\sqrt{2r_j} \label{eq: t and lambda}\\
d_{ij}&=2P_{ij}. \label{eq: d and plucker}
    \end{align}
    
\end{proposition}
\end{tcolorbox}

We defer the proof of Proposition~\ref{prop: scaling} until later. 

\begin{remark}
$\!$The numbers~$P_{ij}$ appearing in \eqref{eq: d and plucker} are the $2\times 2$ minors of the matrix 
\begin{equation*}
\begin{bmatrix}
\cos(\alpha_1) & \cos(\alpha_2) & \cos(\alpha_3) & \cos(\alpha_4) \\
\sin(\alpha_1) & \sin(\alpha_2) & \sin(\alpha_3) & \sin(\alpha_4) 
\end{bmatrix}.    
\vspace{-2pt}
\end{equation*}
This matrix is of a special kind: 
its columns are unit vectors that  lie in the upper~half plane and are ordered counterclockwise. 
Thus, we are not dealing here with the most general form of 
the 3-term Pl\"ucker relation~\eqref{eq: plucker}.
On~the other hand, the general relation~\eqref{eq: plucker} follows from this special case
due to 
\begin{itemize}[leftmargin=.15in]
\item 
the invariance of~\eqref{eq: plucker} under arbitrary rescaling of the columns of the matrix, 
and
\item
the invariance of~\eqref{eq: plucker} under arbitrary permutations of the column indices~$1,2,3,4$ 
(noting that~${P_{ji}=-P_{ij}}$).
\end{itemize}
\end{remark}

\section*{The torus action}

The solutions of the 3-term relation~\eqref{eq:aij} form a quadric hypersurface~$\Sigma\subset\mathbb{C}^6$. 
By Observation~\ref{obs: scaling}, the torus~$(\mathbb{C}^*)^4$ acts on $\Sigma$ by rescaling:
an element 
$(q_1,q_2,q_3,q_4)\in(\mathbb{C}^*)^4$ acts on a 6-tuple~$(a_{ij})\in\Sigma$ 
by sending 
\begin{equation*}
a_{ij}\mapsto q_i q_j a_{ij}. 
\end{equation*}
The action of $(\mathbb{C}^*)^4$ splits the 5-dimensional quadric~$\Sigma$ into 4-dimensional
torus orbits. 

Two points~$(a_{ij}),(b_{ij})\in\Sigma$ will typically lie in different orbits.
In other words, given two 6-tuples~$(a_{ij})$ and~$(b_{ij})$ satisfying~\eqref{eq:aij}--\eqref{eq:bij}, 
there is generally no reason to expect 
the existence of scalars~$q_1,q_2,q_3,q_4$ such that formulas~\eqref{eq:ab-rescaling} hold. 
It is natural to ask: when do such scalars~$q_i$ exist? The answer is given in the following proposition.

\begin{tcolorbox}

\begin{proposition}
\label{prop: rescaling-criterion}
Given two 6-tuples~$(a_{ij})$ and $(b_{ij})$ of nonzero numbers satisfying \eqref{eq:aij}--\eqref{eq:bij}, the following are equivalent:
\begin{itemize}[leftmargin=.15in]
\item 
there exists a quadruple of scalars~$(q_i)$ such that~\eqref{eq:ab-rescaling} holds; 
\item
the numbers~$a_{ij}$ and~$b_{ij}$ satisfy the relation
\begin{equation}
\label{eq:cross-ratios-equal}
\frac{a_{12}a_{34}}{a_{23}a_{14}}=\frac{b_{12}b_{34}}{b_{23}b_{14}} .
\end{equation}
\end{itemize}
\end{proposition}
\end{tcolorbox}

\begin{proof}
It~is evident that~\eqref{eq:ab-rescaling} implies~\eqref{eq:cross-ratios-equal}. It remains to verify the converse implication. 

Suppose that~\eqref{eq:aij}--\eqref{eq:bij} and~\eqref{eq:cross-ratios-equal} hold.  
For~$1\le i<j\le 4$, let us denote
$
c_{ij}=\frac{b_{ij}}{a_{ij}}$.
Equation~\eqref{eq:cross-ratios-equal} implies that
\begin{equation}
\label{eq:c12c34}
c_{12} c_{34} = c_{23} c_{14} . 
\end{equation}
We also have
\begin{equation*}
\frac{a_{13}a_{24}}{a_{23}a_{14}}
\stackrel{\eqref{eq:aij}}{=} \frac{a_{12}a_{34}}{a_{23}a_{14}} + 1
\stackrel{\eqref{eq:cross-ratios-equal}}{=}  \frac{b_{12}b_{34}}{b_{23}b_{14}} + 1
\stackrel{\eqref{eq:bij}}{=} \frac{b_{13}b_{24}}{b_{23}b_{14}} ,
\end{equation*}
implying
\begin{equation}
\label{eq:c23c14}
c_{23} c_{14} = c_{13} c_{24} . 
\end{equation}

We need to find~$q_1,q_2,q_3,q_4$ that satisfy the equations~$q_i q_j = c_{ij}$ (cf.~\eqref{eq:ab-rescaling}). 
To this end, choose~$q_1$ so that 
\begin{equation}
\label{eq:q1}
q_1^2 = \frac{c_{12}c_{13}}{c_{23}} ,
\end{equation}
then set
\begin{equation}
\label{eq:q2q3q4}
q_2 = \frac{c_{12}}{q_1}, \qquad
q_3 = \frac{c_{13}}{q_1}, \qquad
q_4 = \frac{c_{14}}{q_1} .
\end{equation}
We then have~$q_1q_j\!=\!c_{1j}$ (for~$j\!=\!2,3,4$), and moreover 
\begin{align*}
q_2 q_3 & \stackrel{\eqref{eq:q2q3q4}}{=} \frac{c_{12}c_{13}}{q_1^2} \stackrel{\eqref{eq:q1}}{=} c_{23} , \\
q_2 q_4 &\stackrel{\eqref{eq:q2q3q4}}{=} \frac{c_{12}c_{14}}{q_1^2} \stackrel{\eqref{eq:q1}}{=} \frac{c_{14} c_{23}}{c_{13}}
   \stackrel{\eqref{eq:c23c14}}{=} c_{24} , \\
q_3 q_4 &\stackrel{\eqref{eq:q2q3q4}}{=} \frac{c_{13}c_{14}}{q_1^2} \stackrel{\eqref{eq:q1}}{=} \frac{c_{14} c_{23}}{c_{12}}
  \stackrel{\eqref{eq:c12c34}}{=} c_{34}\, ,
\end{align*}
as desired.
\end{proof}

\section*{The cross-ratio}

The reader may recognize the expressions appearing in~\eqref{eq:cross-ratios-equal} as cross-ratios of points on the projective line. 

\begin{tcolorbox}
\begin{definition}
\label{def: cross-ratio}
Let $X_1,X_2,X_3,X_4\in\mathbb{CP}^1$ be distinct points on the complex projective line 
represented by the column vectors $\left[\begin{smallmatrix}x_i \\y_i \end{smallmatrix}\right]$.
    Their \emph{cross-ratio} $[X_1,X_2,X_3,X_4]$ is defined by 
\begin{equation}
\label{eq:cross-ratio}
[X_1,X_2,X_3,X_4] = \frac{P_{12}P_{34}}{P_{23}P_{14}}\,,
\end{equation}
where we use the notation 
\begin{equation*}
P_{ij}=\det\left[\begin{matrix}x_i & x_j\\y_i & y_j\end{matrix}\right]. 
\end{equation*}
This notion can be extended 
to cases where two (but not more) of the $X_i$ coincide.  
\end{definition}
\end{tcolorbox}

The reader is encouraged to verify that Definition~\ref{def: cross-ratio} matches the conventional definition of the cross-ratio.

Definition~\ref{def: cross-ratio} and the last statement of Theorem~\ref{thm: plucker} allow us to restate Proposition~\ref{prop: rescaling-criterion}. 
In the language of cross-ratios, Proposition~\ref{prop: rescaling-criterion} asserts that for two $2\times 4$ matrices $M$ and~$N$, the following are equivalent:
\begin{itemize}[leftmargin=.15in]
\item 
the columns of $N$ can be rescaled so that the Pl\"ucker coordinates match the Pl\"ucker coordinates of~$M$; 
\item
the columns of $M$ and $N$ represent two quadruples of points in $\mathbb{CP}^1$ that have the same cross-ratio. 
\end{itemize}

\section*{Models of the hyperbolic plane}
The proof of Proposition~\ref{prop: scaling} will utilize three different models of the hyperbolic plane.

\subsection*{The upper half plane}
\begin{tcolorbox}
    \begin{definition}
    Our first model of the hyperbolic plane is the \emph{upper half plane}
\begin{equation*}
\mathcal{U}=\{(x,y)\in\R^2\ : \ y>0\}
\end{equation*}
endowed with the Riemannian metric
\begin{equation*}
ds^2=(dx^2+dy^2)/y^2.
\end{equation*}
Identifying~$\R^2$ with the complex plane~$\C$ identifies~$\mathcal{U}$ 
with the complex numbers with positive imaginary part.
The points of the extended real axis~$\R\cup\{\infty\}$ are the \emph{ideal points} of~$\mathcal{U}$. 
Denote by~$\ol{\mathcal{U}}$ the union of~$\mathcal{U}$ and the ideal points.
\end{definition}
\end{tcolorbox}

\begin{figure}[htbp!]
\begin{tikzpicture}[scale=1]
    \draw (-1.2,0) to (1.2,0);
    \draw [black,very thick,domain=0:180] plot ({cos(\x)}, {sin(\x)});
    \node[below] (w1prime) at (-1,0) {$W_1'$};
    \filldraw (-1,0) circle (.04);
    \node[below] (w2prime) at (1,0) {$W_2'$};
    \filldraw (1,0) circle (.04);
    \node[above left] (w1) at ({-sqrt(3)/2},{1/2}) {$W_1$};
    \filldraw ({-sqrt(3)/2},{1/2}) circle (.04);
    \node[above right] (w2) at ({sqrt(2)/2},{sqrt(2)/2}) {$W_2$};
    \filldraw ({sqrt(2)/2},{sqrt(2)/2}) circle (.04);
\end{tikzpicture}  
\vspace{-7pt}
\caption{The ideal points $W_1'$ and $W_2'$ used to compute the hyperbolic distance between $W_1$ and $W_2$}
\vspace{-7pt}
\label{fig:W1W2}
\end{figure}
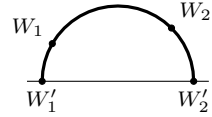

\begin{tcolorbox}
\begin{proposition}[cf.~\cite{penner_decorated_2012}, p.~16]
\label{prop: hyp dist}    
Take two distinct points ${W_1,W_2\in\mathcal{U}}$. 
Let~$W_1',W_2'\in\overline{\mathcal{U}}$ be the ideal endpoints of the geodesic connecting \linebreak[3]
$W_1$ and~$W_2$, see Figure~\ref{fig:W1W2}.  
Then the hyperbolic distance~$\delta_{12}$ between~$W_1$ and~$W_2$ is given by the cross-ratio 
    \[
\delta_{12}=\log[W_1,W_1',W_2,W_2'].\]
\end{proposition}
\end{tcolorbox}

\subsection*{The Poincar{\'e} disk}

\begin{tcolorbox}
    \begin{definition}
    Our second model of the hyperbolic plane is the \emph{Poincar\'e disk}
\begin{equation*}
\DD = \{z\in\C\ : \ |z|=x^2+y^2<1\}
\end{equation*}
endowed with the metric
\begin{equation*}
ds^2=4(dx^2+dy^2)/(1-|z|^2)^2. 
\end{equation*}
Ideal points lie on the unit circle $\bfS$.
\end{definition}
\end{tcolorbox}

To move back and forth between the models~$\mathcal{U}$ and~$\DD$, 
one can apply the mutually inverse isometries known as the \emph{Cayley transforms}:
\begin{align*}
    \mathcal{U} &\to \DD & \DD&\to \mathcal{U}  \\
z&\mapsto \tfrac{z-i}{z+i}&w&\mapsto i\tfrac{1+w}{1-w}.
\end{align*}

\subsection*{The hyperboloid}

\begin{tcolorbox}
\begin{definition}
\emph{Minkowski 3-space} is~$\R^3$ endowed with the indefinite pairing
    \begin{align*}
        \lr{\cdot,\cdot}: \R^3\times\R^3&\to\R\\
        \lr{(x,y,z),(x',y',z')}&\mapsto xx'+yy'-zz'.
    \end{align*}
The \emph{upper sheet of a hyperboloid}
\[
\HH=\left\{u=(x,y,z)\in\R^3 \, : \, \lr{u,u}=-1, \, z>0\right\}
\]
inherits a Riemannian metric from the pairing on Minkowski 3-space, producing a model for the hyperbolic plane.
\end{definition}
\end{tcolorbox}

Identify the closed unit disk~$\DD\cup \bfS$ with the horizontal disk at height zero in Minkowski 3-space.  
Then the following pair of mutually inverse maps give isometries between~$\DD$ and~$\HH$:  
\begin{align*}
    \HH&\to\DD \\
    (x,y,z)&\mapsto \ol{(x,y,z)}=\tfrac{1}{1+z}(x,y) \\
     \DD&\to\HH\\
    (x,y)&\mapsto \tfrac{1}{1-x^2-y^2}(2x,2y,1+x^2+y^2)
\end{align*}

\begin{tcolorbox}
\begin{definition}
The \emph{open positive light cone} $L^+$ is defined by 
\[
L^+=\left\{u=(x,y,z)\in\R^3 \ : \ \lr{u,u}=0, \, z>0\right\}.
\]
\end{definition}
\end{tcolorbox}
For~$u=(x,y,z)\in L^+$ since~$\lr{u,u}=0$, we have
\begin{align*}
    0=\lr{u,u}=\lr{(x,y,z),(x,y,z)}&=x^2+y^2-z^2\, .
\end{align*}
Since~$z>0$, we conclude that~$\sqrt{x^2+y^2}=z$.

The map from~$\HH$ to~$\DD$ is given by central projection from the point~$(0,0,-1)$.
The projection extends to a surjective but not injective natural mapping 
\begin{align*}
    L^+&\to \bfS\\
    (x,y,z)&\mapsto \ol{(x,y,z)}=\tfrac{1}{\sqrt{x^2+y^2}}(x,y,0).
\end{align*}

\section*{Horocycles and lambda lengths}

In order to prove Proposition~\ref{prop: scaling}, we will need to recall some facts about horocycles.

\begin{tcolorbox}
    
\begin{definition}
A \emph{horocycle} $H$ on the hyperbolic plane is a curve of constant curvature such that all perpendicular geodesics are limiting parallel and converge to a single ideal point called the \emph{center} of~$H$.

\end{definition}
\end{tcolorbox}
 In the model of the upper half plane~$\mathcal{U}$, a horocycle~$H$ is either a Euclidean circle tangent to the real axis or a horizontal straight line parallel to the real axis.  
In the former case, the point of tangency is the center of~$H$. 
In the model of Poincar{\'e} disk~$\DD$, a horocycle is a Euclidean circle tangent to~$\bfS$.  
\hbox{The point of tangency is the center of~$H$.}
    
     There is an isomorphism between points of the upper light cone~$L^+$ and the collection of all horocycles in the hyperbolic plane~$\HH$, namely
\begin{align*}
    L^+ &\to \{\text{horocycles in }\HH \}\, ,\\
    u&\mapsto h(u)=\left\{v\in\HH \ : \ \lr{u,v}=-\tfrac{1}{\sqrt{2}}\right\}.
\end{align*}
Let~$u\in L^+$.  Then we can map
\[
u\in L^+ \ \mapsto h(u)\in \mathbb{H} \ \xrightarrow{\text{Cayley transform}} \ \ol{h}(u)\in \mathbb{D}\,.
\]
The center of the horocycle~$\overline{h}(u)$ in~$\DD$  is \[\overline{(x,y,z)}=\tfrac{1}{\sqrt{(x^2+y^2)}}(x,y,0)\]
and the Euclidean radius of~$\overline{h}(u)$ in~$\DD$ is~$\frac{1}{1+z\sqrt{2}}$.

\begin{figure}[htbp!]
    \centering
\begin{tikzpicture}[scale=2.5]
    \draw [black,very thick,domain=200:350] plot ({cos(\x)}, {sin(\x)});
    \draw [black,very thick, dashed,domain=180:200] plot ({cos(\x)}, {sin(\x)});
    \draw [black,thick,dashed,domain=350:370] plot ({cos(\x)}, {sin(\x)});
    \node at (.95,-.4) [right] {$\bfS$};

    \draw ({-sqrt(8)/sqrt(27)},{-2/sqrt(27)}) circle ({1/3});
    \coordinate (A_1) at ({-sqrt(2)/sqrt(3)},{-1/sqrt(3)});
    \filldraw (A_1) circle (.02);
    \node at (-{sqrt(8)/sqrt(27)},{-2/sqrt(27)+1/3}) [above] {$H_1$}; 

    \draw[] ({3/(sqrt(32))},{-3/(sqrt(32))}) circle (1/4);
    \coordinate (A_2) at ({1/sqrt(2)},{-1/sqrt(2)});
            \filldraw (A_2) circle (.02);
    \node at ({3/(sqrt(32))},{-3/(sqrt(32))+1/4}) [above] {$H_2$}; 
    
    \draw (A_1) to [out=37,in=133]  (A_2);
    \draw[very thick, brown] (-.215,-.345) to [out=4,in=160] node[above] {$\delta_{12}$}(0.305,-0.424);

    \node (gamma) at (-.5,-.5) {$\gamma$};
\end{tikzpicture}    
\vspace{-3pt}
\caption{The signed hyperbolic distance between two horocycles
}
\vspace{-5pt}
    \label{fig: lambda length}
\end{figure}
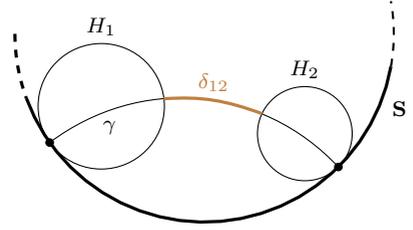

\begin{tcolorbox}
    
\begin{definition}
Let~$H_1,H_2\subset\DD$ be  \hbox{horocycles} with distinct centers.  
Let~$\gamma$ be \hbox{the geodesic that} connects these centers.  Let~$\delta_{12}$ be the~signed hyperbolic distance along~$\gamma$ between the points $H_1\cap\gamma$ and~$H_2\cap\gamma$, where the sign of~$\delta_{12}$ is positive if and only if~$H_1$ and $H_2$ are disjoint, cf.\ Figure~\ref{fig: lambda length}. \linebreak[3]
The \emph{lambda length} between $H_1$ and~$H_2$ is then defined~by
    \[
    \lambda(H_1,H_2)=\lambda_{12}=\sqrt{e^\delta}.
    \]

\end{definition}
\end{tcolorbox}

    Alternatively, let
\begin{equation}
\label{eq:u1u2}
u_1=(x_1,y_1,z_1), \, u_2=(x_2,y_2,z_2)\in L^+
\end{equation}
be horocycles that do not lie on a common ray in~$L^+$. 
The lambda length~$\lambda_{12}=\lambda(u_1,u_2)$ is given by 
\[
\lambda_{12}=\sqrt{-\lr{u_1,u_2}} = \sqrt{-x_1x_2-y_1y_2+z_1z_2}\,.
\]

\section*{Proof of Proposition~\ref{prop: scaling}
}
\label{sec: proofs}

\begin{proof}[Proof of~\eqref{eq: d and t}.]
This equation can be deduced from the law of cosines, see, 
e.g., \cite{zacharias_caseysche_1942}, \cite[p.~297]{fukagawaSacredMathematics2008}.
\end{proof}

\begin{proof}[Proof of~\eqref{eq: t and lambda}.]

Consider horocycles $u_1, u_2\in L^+$, as in~\eqref{eq:u1u2}, 
that do not lie on a common ray in~$L^+$.  
Viewing them as  horocycles~$H_1=\overline{h}(u_1),H_2=\ol{h}(u_2)$ in~$\DD$, we see that they have hyperbolic centers 
\[
A_1=\biggl(\tfrac{x_1}{\sqrt{x_1^2+y_1^2}},\tfrac{y_1}{\sqrt{x_1^2+y_1^2}}\biggr), 
A_2=\biggl(\tfrac{x_2}{\sqrt{x_2^2+y_2^2}},\tfrac{y_2}{\sqrt{x_2^2+y_2^2}}\biggr).
\]
As Euclidean circles, $H_1$ and $H_2$ have radii 
\[
r_1=\tfrac{1}{1+z_1\sqrt{2}}, \quad r_2=\tfrac{1}{1+z_2\sqrt{2}}.
\]
Then
\begin{align*}
\noalign{$(1-r_1)(1-r_2)$\hfill}
    &=1-\tfrac{1}{1+z_1\sqrt{2}}-\tfrac{1}{1+z_2\sqrt{2}}+\tfrac{1}{(1+z_1\sqrt{2})(1+z_2\sqrt{2})}\\
    &=\tfrac{2z_1z_2}{(1+z_1\sqrt{2})(1+z_2\sqrt{2})} .
\end{align*}
The Euclidean distance between~$A_1$ and~$A_2$ is
\begin{align*}
    d_{12}&=\biggl(
    \biggl(\hspace{-2pt}\tfrac{x_1}{\sqrt{x_1^2+y_1^2}}-\tfrac{x_2}{\sqrt{x_2^2+y_2^2}}\hspace{-2pt}\biggr)^{\!\!2}
    \!\!+\biggl(\hspace{-2pt}\tfrac{y_1}{\sqrt{x_1^2+y_1^2}}-\tfrac{y_2}{\sqrt{x_2^2+y_2^2}}\hspace{-2pt}\biggr)^{\!\!2}
    \biggr)^{\!\!1/2}\\
    &=\biggl(
    2\biggl(1-\tfrac{x_1x_2+y_1y_2}{\sqrt{x_1^2+y_1^2}\sqrt{x_2^2+y_2^2}}\biggr)\biggr)^{\!\!1/2} \\
    &=\Bigl(2
    \left(1-\tfrac{x_1x_2+y_1y_2}{z_1z_2}\right)
    \Bigr)^{\!1/2} \\
    &=\Bigl(2
    \left(\tfrac{z_1z_2-x_1x_2-y_1y_2}{z_1z_2}\right)
    \Bigr)^{\!1/2}\, .
\end{align*}
By equation~\eqref{eq: d and t}, the length~$t_{12}$ of the bitangent is
\begin{align*}
    t_{12}
    &=\Bigl(2 \left(\tfrac{z_1z_2-x_1x_2-y_1y_2}{z_1z_2}\right)\left(\tfrac{2z_1z_2}{(1+z_1\sqrt{2})(1+z_2\sqrt{2})}\right)\Bigr)^{\!1/2}\\
    &=2(-x_1x_2-y_1y_2+z_1z_2)^{1/2} ((1+z_1\sqrt{2}) (1+z_2\sqrt{2}))^{\!-1/2} \\
    &=2\lambda_{12}\cdot ((1+z_1\sqrt{2}) (1+z_2\sqrt{2}))^{\!-1/2} \\
    &=\lambda_{12}\sqrt{2r_1}\sqrt{2r_2}.\qedhere
\end{align*}
\end{proof}
\begin{figure*}[htbp!]
    \centering
\begin{tikzpicture}[scale=4]
    \draw [very thick,domain=200:350] plot ({cos(\x)}, {sin(\x)});
    \draw [very thick, dashed,domain=180:200] plot ({cos(\x)}, {sin(\x)});
    \draw [thick,dashed,domain=350:370] plot ({cos(\x)}, {sin(\x)});

    \filldraw[] (0,0) circle  (.02);
    \coordinate (o) at (0,0);
    \node at (o) [above] {$O$};

    \coordinate (A_i) at ({-sqrt(2)/sqrt(3)},{-1/sqrt(3)});
    \node at (A_i) [below left] {$A_i=(\cos2\alpha_i,\sin2\alpha_i)$};

    \coordinate (A_j) at ({1/sqrt(2)},{-1/sqrt(2)});
    \node at (A_j) [below right] {$A_j=(\cos 2\alpha_j,\sin 2\alpha_j)$};    
    \coordinate (mid) at ({(sqrt(3)-2)/(2*sqrt(6))},{(-sqrt(3)-sqrt(2))/(2*sqrt(6))});

    \draw[] (mid) -- 
    ({-sqrt(2)/sqrt(3) + (1/sqrt(2)+sqrt(2)/sqrt(3))/20*9},{-1/sqrt(3)+(-1/sqrt(2)+1/sqrt(3))/20*9})--
    ({-sqrt(2)/sqrt(3) + (1/sqrt(2)+sqrt(2)/sqrt(3))/20*9-(-1/sqrt(2)+1/sqrt(3))/20*1},{-1/sqrt(3)+(-1/sqrt(2)+1/sqrt(3))/20*9+(1/sqrt(2)+sqrt(2)/sqrt(3))/20*1}) -- ({-sqrt(2)/sqrt(3) + (1/sqrt(2)+sqrt(2)/sqrt(3))/20*9-(-1/sqrt(2)+1/sqrt(3))/20*1+(1/sqrt(2)+sqrt(2)/sqrt(3))/20*1},{-1/sqrt(3)+(-1/sqrt(2)+1/sqrt(3))/20*9+(1/sqrt(2)+sqrt(2)/sqrt(3))/20*1+(-1/sqrt(2)+1/sqrt(3))/20*1}) -- (mid);

    \draw (mid) to node[right, near start] {$ \cos(\alpha_j-\alpha_i)$} (o);

    \draw[very thick, magenta] (A_i) to node[below, pos=.6] {$ \frac{d_{ij}}{2}=\sin(\alpha_j-\alpha_i)$} (mid);
    \draw (mid) to (A_j);

    \draw (A_i) to node[above left] {$1$} (o);

    \draw (A_j) to node[above right] {$1$} (o);

\begin{scope}
    \clip (o) to (A_i) to (mid) to (o);
    \draw (o) circle (.1);
    \node at (o) [left, rotate=60,inner sep=15pt] {$ \alpha_j-\alpha_i$};
  \end{scope}

\end{tikzpicture}    
\captionsetup{width=9cm}

\vspace{-3pt}
\caption{Measurements used to calculate~$d_{ij}$}
\vspace{-3pt}
    \label{fig: plucker setup}
\end{figure*}
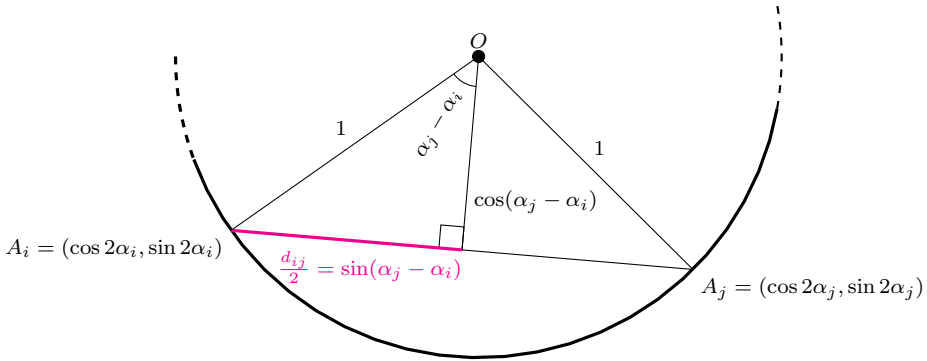

\begin{proof}[Proof of~\eqref{eq: d and plucker}.]
Recall that $A_i=(\cos2\alpha_i,\sin2\alpha_i)$ and $\alpha_i<\alpha_j$ for~$i<j$. 
We then compute (see Figure~\ref{fig: plucker setup}):
\begin{align*}
    d_{ij}&=2\sin(\alpha_j-\alpha_i)\\
    &=2(\sin \alpha_j\cos \alpha_i-\cos \alpha_j\sin \alpha_i)\\
    &=2\left[
    \begin{matrix}
        \cos \alpha_i & \cos \alpha_j \\
        \sin \alpha_i & \sin \alpha_j
    \end{matrix}
    \right
    ]\\
    &=2P_{ij}\, .\qedhere
\end{align*}

\end{proof}

\section*{Generalizations and historical notes}
\label{sec: further results} 

Ptolemy's theorem (Theorem~\ref{thm: ptolemy}) can be generalized in multiple directions. 
In 1872, Darboux \cite{darbouxRelationsEntreGroupes1872} shifted the setting to spherical geometry, establishing a relation among pairwise spherical distances between  four concyclic points on a sphere; similar results were obtained by Frobenius~\cite{frobenius_anwendungen_1875}. In 1912, Kubota~\cite{kubota_extended_1912} found an analogous relation in a hyperbolic setting. Indeed, Observation~\ref{obs: scaling} can be used to extend Ptolemy's relation to both spherical and hyperbolic geometry, see Guo and Sonmëz \cite{guoCyclicPolygonsClassical2011}, the discussion in \cite[p.43]{penner_decorated_2012}, and the pictorial explanation on Stothers' remarkable website~\cite{stothersPtolemysTheoremHyperbolic}.

Various authors looked at the converse of Theorem~\ref{thm: ptolemy}: 
four points only lie on a circle if the 3-term relation holds among their pairwise distances. Haantjes~\cite{haantjesCharacteristicLocalProperty1951} gave both spherical and hyperbolic versions of this result; 
this was later completed by Valentine~\cite{valentine_analogue_1970,valentine_analogue_1970-1}.

One could  generalize Theorem~\ref{thm: ptolemy} by considering larger collections of concyclic points.  
A result for hexagons was obtained in the late 19th century by Fuhrmann~\cite{fuhrmann_synthetische_1890}. 
This was further extended by Gregorac~\cite{gregorac_feuerbachs_1996}.

Casey's relation holds in greater generality than is stated in Theorem~\ref{thm: casey},
for any configuration of four circles that are tangent---either internally or externally---to a fifth circle,  
as long as one is careful about their choice of bitangents. 
For a comprehensive summary, see~\cite{gueron_two_2002}.  

Various authors have found versions of Casey's theorem in spherical and hyperbolic geometry \cite{abrosimovCaseysTheoremHyperbolic2015,kostin_interpretation_2016} and in higher dimensions \cite{abrosimov_generalizations_2018,maehara_bipartite_2019}. 
Kostin~\cite{kostin_analogs_2024} extended Fuhrmann's theorem to the hyperbolic plane.

Mathews (resp., Mathews and Varsha) extended Penner’s lambda lengths to hyperbolic space of dimension~3 (resp.,~4) using complex numbers
(resp., quaternions), and demonstrated that these extensions satisfy
equations similar to Theorem~\ref{thm: penner}, see \cite{mathews_spinors_2024,mathews_spinors_2025,mathews_quaternionic_2025}.

An uninitiated reader may be interested in learning more about Shiraishi Nagatada (or Shiraishi Ch\=och\=u)  and the practice of hanging \emph{sangaku}, wooden tablets with geometric figures, on the walls of Shinto shrines and Buddhist temples in Japan.  This practice was widespread from the 17th through 20th century, coinciding largely with the Edo period when Japan was closed to trade and intellectual exchange with Europe. 
An excellent source to learn more is the survey~\cite{fukagawaSacredMathematics2008}, see also~\cite{hosking_sangaku_2016}. 

Shiraishi's work \emph{Shamei Sanpu} \cite{shiraishi_nagatada_shamei_1826} contains a large collection of sangaku.
Shiraishi was a \emph{samurai}, a member of nobility. At the time, this elite strata of society included mathematicians who also taught students \cite{smith_history_1914}.

Cross-ratios appear in Book 7 of Pappus' \emph{Collection}, which drew heavily on now lost works of Euclid \cite{jones_book_1986}.
Like Ptolemy, both Euclid and Pappus lived and worked in Alexandria, around 300 BCE and 300 CE respectively.
The modern notion of a cross-ratio crystallized during the first half of the 19th century in the work of Lazare  Carnot, Jean-Victor Poncelet,  Michel Chasles, and August M{\"o}bius. 
The cross-ratio plays a central role in the foundations of projective geometry as developed by Karl von Staudt, Arthur Cayley, and Felix Klein. 
Klein also used the cross-ratio to define the metric on hyperbolic space, see, e.g., \cite[Proposition~\ref{prop: hyp dist}]{laptev_mathematics_1996}.

Theorems~\ref{thm: ptolemy}-\ref{thm: plucker} underlie areas of active current research, such as the study of 
cluster algebras and frieze patterns.  
Frieze patterns are arrays of numbers that satisfy a  certain local algebraic relation. 
This relation can be interpreted as an instance of any of the four incarnations
\eqref{eq: ptolemy}--\eqref{eq: plucker} of the 3-term relation $AB+CD=EF$.
When H.~S.~M.~Coxeter \cite{coxeterFriezePatterns1971} introduced frieze patterns, 
he observed their connection to cross-ratios and to the \emph{Pentagramma mirificum} of Gauss and Napier.  

As shown by Coxeter and Conway \cite{coxeterFriezePatterns1971,conway_triangulated_1973a}, 
frieze patterns can be studied from a combinatorial perspective that utilizes triangulations of a polygon. 
This construction leads to at least two geometric interpretations of friezes. 
First, Penner's Theorem~\ref{thm: penner} can be used to realize frieze patterns 
(and more generally, cluster algebras associated with bordered oriented surfaces) 
via lambda lengths for triangulations of ideal hyperbolic polygons \cite{fomin_cluster_2018,felikson_ptolemy_2023} 
(resp., triangulations of the surface) decorated by horocycles.  
Second, equation~\eqref{eq: t and lambda} and Casey's formula allow us to realize Coxeter-Conway friezes via bitangent distances between four circles tangent to a fifth circle. 

\newpage

\section*{Acknowledgements}
I am grateful to Michael Shapiro, who asked whether Theorems~\ref{thm: casey} and~\ref{thm: penner} can be obtained from each other, which began the exploration that led to this note. I thank Sergey Fomin for sharing his observation that the diamond rule in a frieze pattern can be geometrically interpreted as Casey's formula, and many mathematical and editorial discussions that helped organize my thinking.

\begin{flushleft}
\bibliographystyle{myamsalpha}
\bibliography{herfriezereferences.bib}
\RaggedRight
\end{flushleft}
\end{document}